\documentclass[12pt]{article}
\usepackage{graphicx} 
\usepackage{amsmath,amssymb,amsthm}
\usepackage{hyperref}
\usepackage{tikz}
\usetikzlibrary{patterns}
\usepackage{authblk}

\title{quivers}
\author[1]{Tam\'as Terpai}
\affil[1]{Department of Analysis, Institute of Mathematics, ELTE Eötvös Loránd University, Budapest, Hungary}

\newcommand{\rk}{{\operatorname{rank}}}
\newcommand{\Z}{{\mathbb{Z}}}
\newcommand{\ZN}{{\mathbb{Z}_{\geq 0}}}

\newtheorem{lemma}{Lemma}
\newtheorem{thm}{Theorem}
\newtheorem*{remark}{Remark}

\begin{document}

\title{A note on numerical symmetry of the product-zero variety}
\maketitle

\begin{abstract}
We provide an elementary proof of a symmetric property of the top dimensional components of the quiver subvariety $A_n\dots A_1=0$.
\end{abstract}

\section{Introduction}

Given a dimension vector $\vec d = (d_0,d_1,\dots,d_n)$ one can consider a type $A$ quiver $V_0 \to V_1 \to \dots \to V_n$ with complex vector spaces $V_j$ of dimension $\dim V_j=d_j,0 \leq j \leq n$, equipped with the natural action of $G=G_{\vec d}=GL(V_0) \times \dots \times GL(V_n)$ on $H=H_{\vec d}=Hom(V_0,V_1) \times \dots \times Hom(V_{n-1},V_n)$, and define the product-zero variety
\[
\Sigma=\Sigma_{\vec d} =\{ (A_1,\dots,A_n) \in H_{\vec d} \mid A_n\dots A_1=0\}.
\]
This variety is studied in detail in \cite{topDim} and \cite{symm} and one surprising result is that the number and dimension of the top dimensional components of $\Sigma_{\vec d}$ are the same for those $\vec d$ whose entries form permutations of one another. The proofs involve highly abstract tools, and the main purpose of this note is to provide an elementary proof.

In Section \ref{sec:trunk} we reformulate the problem in a purely combinatorial way and reduce it to a technical statement that is proved in Section \ref{sec:comm}. 

\section{Reduction to combinatorics}\label{sec:trunk}

To identify the components of $\Sigma$, we first decompose it into finitely many subvarieties. For a fixed integer \emph{rank sequence} $\vec r=(r_0,r_1,\dots,r_n)$, $d_0=r_0 \geq r_1 \geq \dots \geq r_n \geq 0$, define
\[
H_{\vec d}(\vec r)= \{ (A_1,\dots,A_n) \in H_{\vec d} \mid \forall j: \rk A_j\dots A_1= r_j \}.
\]
We obviously have $H_{\vec d}=\underset{\vec r}{\cup} H_{\vec d}(\vec r)$ and $\Sigma_{\vec d}=\underset{\vec r: r_n=0}{\cup} H_{\vec d}(\vec r)$. Note that all $H_{\vec d}(\vec r)$, while not necessarily orbits of $G_{\vec d}$, are nevertheless irreducible. Indeed, using induction, for $n=1$ they are in fact the orbits of $G$ and hence irreducible, and for $\vec d=(d_0,\dots,d_n)$, $\vec{d'}=(d_0,\dots,d_{n-1})$, $\vec r=(r_0,\dots,r_n)$, $\vec{r'}=(r_0,\dots,r_{n-1})$ we have a natural forgetful map $\rho: H_{\vec d}(\vec r) \to H_{\vec{d'}}(\vec{r'})$. This $\rho$ is surjective whenever $H_{\vec d}(\vec r)\neq \emptyset$ (that is, $r_n \leq d_n$ and $H_{\vec{d'}}(\vec{r'})\neq \emptyset$) and in this case it is a fiber bundle with irreducible fiber
\[
\{ A_n \in Hom(V_{n-1},V_n) \mid rk (A_n|_{\operatorname{im} A_{n-1}\dots A_1})=r_n\} \cong
\]
\[
\cong \{B \in Hom(\mathbb C^{r_{n-1}},\mathbb C^{d_n}) \mid \rk B=r_n \} \times Hom(\mathbb C^{d_{n-1}-r_{n-1}},\mathbb C^{d_n}).
\]

This subdivision shows that the components of $\Sigma_{\vec d}$ can only be closures of those $H_{\vec d}(\vec r)$ that do not lie in the closure of some other stratum. In particular, the top dimensional components of $\Sigma_{\vec d}$ are in $1$-to-$1$ correspondence with the top dimensional $H_{\vec d}(\vec r)$, and from now on we will try to count these. It is certainly possible to track the lower dimensional components as well in a case-by-case analysis, but the condition that an $H_{\vec d}(\vec r)$ is not in the closure of another one -- whenever $r_k<r_{k-1}$, we also have $d_{k+1}-r_{k+1}\geq d_k-r_k$ -- makes a general approach complicated enough that we do not have practical formulas for this setting.

To select the top dimensional $H_{\vec d}(\vec r)$, we first calculate the codimension of $H_{\vec d}(\vec r)$ in $H_{\vec d}$ using the sequence of forgetting maps $\rho$ from before:
\begin{align*}
\delta(\vec d,\vec r) \overset{\text{def}}{=} \operatorname{codim} H_{\vec d}(\vec r) &= \sum_{j=1}^n \operatorname{codim} \{B \in Hom(\mathbb C^{r_{j-1}},\mathbb C^{d_j}) \mid \rk B=r_j \} =\\
&= \sum_{j=1}^n (r_{j-1}-r_j)(d_j-r_j).
\end{align*}
If we again use $\vec{d'}=(d_0,\dots,d_{n-1})$ and $\vec{r'}=(r_0,\dots,r_{n-1})$ to denote the starting slices of $\vec d$ and $\vec r$, respectively, then we obtain a recursive form
\begin{equation}\tag{A}\label{eq:rec}
\delta(\vec d,\vec r) = \delta(\vec{d'},\vec{r'})+(r_{n-1}-r_n)(d_n-r_n).
\end{equation}
We now define a generating function
\[
F_{\vec d}(x,y)=\sum_{\vec r} x^{r_n} y^{\delta(\vec d,\vec r)}
\]
that we consider as a polynomial in $\ZN[y][x]$. The original question thus becomes finding the minimal degree monomial in the ``constant'' ($x$-degree $0$) term of $F_{\vec d}$, its coefficient will be the number of top dimensional components and its $y$-degree will be the codimension in $H_{\vec d}$ of these components. The key statement that we proceed to prove now is the following:

\begin{lemma}\label{lemma:symm}
The leading (minimal $y$-degree) terms of the coefficients of $F_{\vec d}$ are invariant under permuting the entries of $\vec d$.
\end{lemma}
\noindent As an immediate corollary we obtain the symmetry result of \cite{topDim} and \cite{symm}:
\begin{thm}
If the entries of $\vec{d_1}$ form a permutation of the entries of $\vec{d_2}$, then the dimension of $\Sigma_{\vec{d_1}}$ is the same as that of $\Sigma_{\vec{d_2}}$ and the number of top dimensional components of $\Sigma_{\vec{d_1}}$ is the same as that of $\Sigma_{\vec{d_2}}$.
\end{thm}

To make the arguments easier to follow, we will use the notation $f\sim g$ to mean that the polynomials $f,g \in \Z[y][x]$ have the same leading terms in their coefficients (or, in other words, that they agree up to higher order terms). Using this notation, the statement of Lemma \ref{lemma:symm} can be reformulated as ``for any permutation $\sigma$ of the index set $\{0,1,\dots,n\}$ we have $F_{(d_0,\dots,d_n)} \sim F_{(d_{\sigma(0)},\dots,d_{\sigma(n)})}$''.

To prove Lemma \ref{lemma:symm}, we split the defining sum of $F_{\vec d}$ according to $r_n$ and use \eqref{eq:rec} to get
\begin{align*}
F_{\vec d}(x,y)&=\sum_{r_n,\vec{r'}} x^{r_n}y^{\delta(\vec{d'},\vec{r'})+(r_{n-1}-r_n)(d_n-r_n)}=\\
&=\sum_{\vec{r'}} y^{\delta(\vec{d'},\vec{r'})}\sum_{r=0}^{\min\{r_{n-1},d_n\}} x^ry^{(r_{n-1}-r)(d_n-r)}.
\end{align*}
This expression can be obtained from the defining sum of $F_{\vec{d'}}$ by $y$-linearly replacing each $x^s$ -- that come from the $\vec{r'}$ having $r_{n-1}=s$ -- by the sum
\[
\Phi_{d_n}(x^s) \overset{\text{def}}{=} \sum_{r=0}^{\min\{s,d_n\}} x^ry^{(s-r)(d_n-r)}.
\]
That is, if we extend this definition of $\Phi_{d_n}$ to a $y$-linear map on $\Z[y][x]$, then
\[
F_{\vec d}(x,y)=\Phi_{d_n}(F_{\vec{d'}}(x,y))
\]
and continuing this reduction we arrive at the expression
\[
F_{\vec d}(x,y)=\Phi_{d_n}(\Phi_{d_{n-1}}(\dots \Phi_{d_1}(F_{(d_0)}(x,y))\dots))=\Phi_{d_n}(\Phi_{d_{n-1}}(\dots \Phi_{d_1}(x^{d_0})\dots)).
\]
Here $d_0$ seems to play a different role than the rest of the $d_j$'s, but we can make it similar to the rest by replacing $x^{d_0}$ by $\Phi_{d_0}(x^N)$ for a sufficiently large $N$. This change adds new terms, but all of those terms will carry a coefficient of $y$-degree at least $N-d_0+1$. Applying further $\Phi_{d_j}$ cannot decrease this degree, so the end result will also get only extra terms with $y$-degree at least $N-d_0+1$, which can be made arbitrarily large and not change the leading terms of $F_{\vec d}(x,y)$:
\[
F_{\vec d}(x,y)\sim\Phi_{d_n}(\Phi_{d_{n-1}}(\dots \Phi_{d_0}(x^N)\dots)).
\]
A similar argument shows that $y$-linearity of $\Phi_d$ implies that it is compatible with the relation $\sim$ in the sense that if $f\sim g$ (and $f,g\in \ZN[y][x]$), then $\Phi_d(f)\sim\Phi_d(g)$. Indeed, if a leading term $x^ay^b$ in a polynomial $f$ produces a term $x^\alpha y^\beta$ in $\Phi_d(f)$, then adding a higher order term $x^ay^{b+\Delta}$ modifies the result by $x^\alpha y^{\beta+\Delta}$, which is superceded by $x^\alpha y^\beta$ (and there is no cancellation to get around this).

This means that to finish the proof of Lemma \ref{lemma:symm} it is enough to show that for any polynomial $f\in \ZN[y][x]$ and indices $a,b$ we have
\[
\Phi_a(\Phi_b(f))\sim \Phi_b(\Phi_a(f)).
\]
This is the commutation relation that we prove in the following section.

\section{Commutation (up to higher order terms)}\label{sec:comm}

When written in the natural basis $1,x,x^2,\dots$, the transformation $\Phi_L$ for an arbitrary $L$ has the following matrix:
\[
\left[\Phi_L\right]=
\begin{pmatrix}
1 & y^L & y^{2L} & y^{3L} & \cdots & y^{(L-1)L} & y^{L^2} & y^{(L+1)L} & \cdots \\
0 & 1 & y^{L-1} & y^{2L-2} & \cdots & y^{(L-2)(L-1)} & y^{(L-1)^2} & y^{L(L-1)} & \cdots \\
0&0&1 & y^{L-2} & \cdots & y^{(L-3)(L-2)} & y^{(L-2)^2} & y^{(L-1)(L-2)} & \cdots \\
0&0&0&1  & \cdots & y^{(L-4)(L-3)} & y^{(L-3)^2} & y^{(L-2)(L-3)} & \cdots \\
\vdots &\vdots&\vdots&\vdots & \ddots & \vdots &\vdots&\vdots\\
0 & 0 & 0 & 0 & \cdots & 1 & y & y^2 & \cdots \\
0 & 0 & 0 & 0 & \cdots & 0 & 1 & 1 & \cdots \\
0 & 0 & 0 & 0 & \cdots & 0 & 0 & 0 & \cdots \\
\vdots & \vdots & \vdots & \vdots &  & \vdots & \vdots & \vdots & \ddots 
\end{pmatrix}
\]
Its entry at $(i+1,j+1)$, $a^{(L)}_{i,j}$ is equal to $y^{(L-i)(j-i)}$ whenever $(0\leq)i\leq j$ and $i\leq L$ (and is $0$ otherwise).

To make the patterns of the obtained matrices more similar to one another, we first change the basis $1,x,x^2,\dots,x^n,\dots$ to $1,x,x^2y,\dots,x^ny^{\binom{n}{2}},\dots$. This causes the entry $(i+1,j+1)$ in $\left[\Phi_L\right]$ to get multiplied by $y^{\binom{j}{2}}y^{-\binom{i}{2}}$, thus it becomes
\begin{align*}
b^{(L)}_{i,j}&=y^{(L-i)(j-i)}y^{\frac{j(j-1)}{2}-\frac{i(i-1)}{2}}=y^{\frac{j^2-i^2-j+i+2(L-i)(j-i)}{2}}=\\&=y^{\frac{(j-i)(2(L-i)+i+j-1)}{2}}=y^{\frac{(j-i)(2L-1+(j-i))}{2}}
\end{align*}
(whenever not $0$ due to $i>j$ or $i>L$). In particular, the diagonals parallel to the main diagonal all carry the same value in them until it becomes $0$ at rows $L+2$ (corresponding to $x^{L+1}$) and above.

Next, we multiply the new matrix $((b^{(L)}_{i,j}))_{i,j\geq 0}$ by $y^{\binom{L}{2}}$. The new entry at $(i+1,j+1)$ becomes thus
\[
c^{(L)}_{i,j}=y^{\frac{(j-i)^2+(2L-1)(j-i)+L^2-L}{2}}=y^{\frac{(j-i+L)(j-i+L-1)}{2}},
\]
so all the nonzero rows of the resulting matrix
\[
C_L = ((c^{(L)}_{i,j}))_{i,j\geq 0}
\]
are end slices of the same sequence $(1,1,y,y^3,\dots,y^{\binom{n}{2}},\dots)$, starting from the term $n=L$ and shifted to take $y^{\binom{L}{2}}$ onto the main diagonal. What we did so far only multiplied the investigated maps by scalars, so the $\Phi_L$ commute up to higher order terms if and only if the $C_L$ commute up to higher order terms.

We can now reinterpret $C_L$ as a truncated polynomial multiplication operator. Identify the row vector $(p_0(y),\dots,p_k(y)) \in \ZN[y,y^{-1}]^{k+1}$ with the polynomial $p(y,t)=p_0(y)+p_1(y)t+\dots+p_k(y)t^k$ in a new formal variable $t$. Then \emph{right} multiplication by $C_L$ takes $t^j$ to $t^j\sum_{s=0}^\infty y^{\binom{L+s}{2}}t^s$ for $j\leq L$ and to $0$ for $j>L$, so it is taking the degree at most $L$ part of the input polynomial and multiplies it with a fixed polynomial.

To simplify presentation, for any formal power series (or polynomial) $f(x,t)=\sum_{j=0}^\infty f_j(y)t^j \in \Z[y,y^{-1}][[t]]$ let the truncation of $f$ to the interval $[A,B]$ be
\[
f|_A^B =\sum_{j=A}^B f_j(y) t^j,
\]
with the additional convention that $f|_A = f|_A^\infty$ and $f|^B=f|_0^B$. With this notation, the multiplicand $\sum_{s=0}^\infty y^{\binom{L+s}{2}}t^s$ can be rewritten as $t^{-L}q|_L$ for the fixed polynomial
\[
q(y,t)=\sum_{j=0}^\infty y^{\binom{j}{2}}t^j
\]
and right multiplication by $C_L$ takes $p$ to $p|^L\cdot t^{-L}q|_L$.

Returning to the commutation relation, the $j+1^{st}$ row of the product $C_LC_M$ is the image of the $j+1^{st}$ basis vector and hence corresponds to the polynomial
\[
p(j,L,M) \overset{\text{def}}{=} (t^j|^L \cdot t^{-L}q|_L)|^M \cdot t^{-M}q|_M.
\]
If we show that $p(j,L,M)\sim p(j,M,L)$ for all $j$, $L$ and $M$, then the matrix products $C_LC_M$ and $C_MC_L$ agree up to higher order terms and Lemma \ref{lemma:symm} is done.

It is easy to see that if $j>\min\{L,M\}$, then $p(j,L,M)=0$ due to one of the $|^L$ and $|^M$, so we assume now that $j \leq \min \{ L,M \}$. In this case
\begin{align*}
(t^j|^L\cdot t^{-L}q|_L)|^M t^{-M} q|_M&=(t^j\cdot t^{-L}q|_L)|^M t^{-M} q|_M=\\
&=t^j\cdot (t^{-L}q|_L)|^{M-j} t^{-M} q|_M=\\
&=t^jt^{-L}\cdot q|_L^{M+L-j} t^{-M} q|_M=\\
&=t^{j-L-M}\cdot q|_L^{M+L-j} \cdot q|_M.
\end{align*}
The shift by $j-L-M$ is the same when $L$ and $M$ are switched, so it is enough to show that
\begin{equation}\tag{B}\label{eq:qSymm}
q|_L^{M+L-j} \cdot q|_M \sim q|_M^{M+L-j} \cdot q|_L.
\end{equation}
To get the coefficient of degree $D$ of the product on the left-hand side, say, we need to take the points of the domain $[L,M+L-j] \times [M,\infty) \subset \Z^2\langle u,v\rangle$ that lie on the line $u+v=D$ and add together the products $y^{\binom{u}{2}}y^{\binom{v}{2}}=y^{\binom{u}{2}+\binom{v}{2}}$ of the corresponding coefficients of $q$ for all of these points $(u,v)$. Since $u\mapsto \binom{u}{2}$ is strictly convex, the least exponent will occur at the point(s) where $|u-v|$ is minimal, closest to the diagonal $u=v$.

Let us now compare the two domains for $q|_L^{M+L-j} \cdot q|_M$ and $q|_M^{M+L-j} \cdot q|_L$, respectively, and switch the coordinates in the latter one (this doesn't change the weights $y^{\binom{u}{2}+\binom{v}{2}}$):
\[
U=[L,M+L-j] \times [M,\infty) \qquad\mbox{and}\qquad V=[L,\infty) \times [M,M+L-j].
\]
They are both (one-side infinite) rectangles, with bottom-leftmost vertex at $(L,M)$. The intersection of the infinite rays of the boundary opposite to this vertex lies at $(M+L-j,M+L-j)$. This means that up to $u+v=2(M+L-j)$ the common part is closer to the diagonal $u-v=0$ than the parts of the symmetric difference $U\Delta V$, and at higher sums $u+v$, those endpoints of the sections of $U$ and $V$ that are closer to the diagonal are symmetric with respect to the diagonal, so the optimal value of the exponent of $y$ is the same.
\begin{figure}
\centering
\begin{tikzpicture}
  \draw[very thin,gray] (-0.5,-0.5) grid (9.5,9.5);
  \draw[->] (0,-0.5) --
    (0,3) node[fill=white] {$M$} --
    (0,7) node[fill=white] {$M+L-j$} -- (0,9.5);
  \draw[->] (-0.5,0) --
    (5,0) node[fill=white,rotate=90] {$L$} --
    (7,0) node[fill=white,rotate=90] {$M+L-j$} -- (9.5,0);
  \draw[fill,pattern=north east lines] (5,9.5) -- (5,3) -- (7,3) -- (7,9.5);
  \draw (6,8) node {$U$};
  \draw[fill,pattern=north west lines] (9.5,7) -- (5,7) -- (5,3) -- (9.5,3);
  \draw (8,5) node {$V$};
  \draw[thick, color=red] (0,0) -- (9.5,9.5)
    node[near start, above, rotate=45] {$u=v$};
\end{tikzpicture}
\caption{Exponents of $y$ in the two products of \eqref{eq:qSymm}.}
\end{figure}
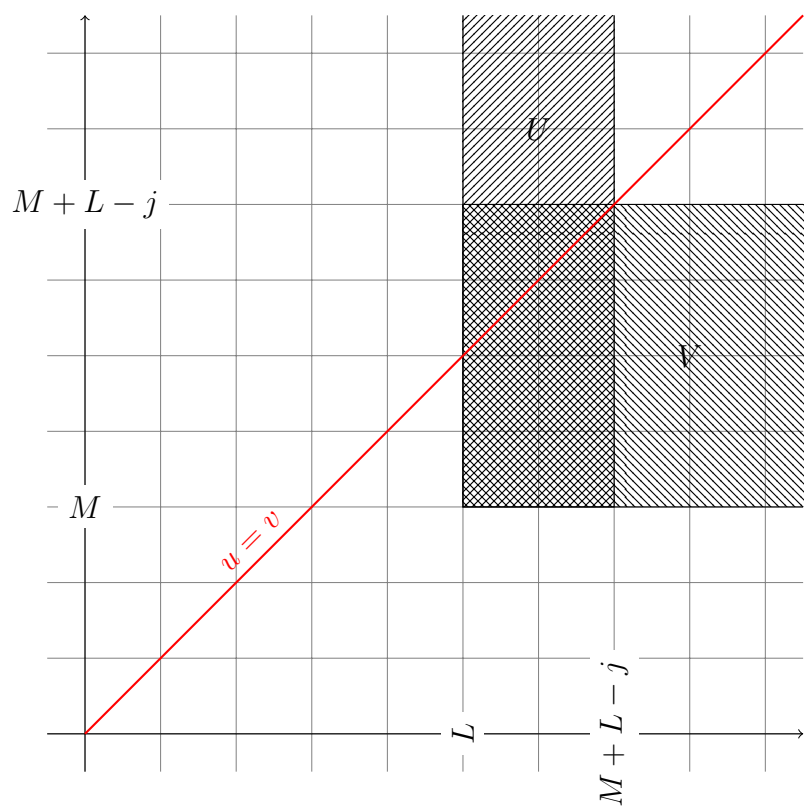

The proof of Lemma \ref{lemma:symm} is finished.

\begin{remark}
While this argument does not provide any deep reason as to why the symmetry holds, it has the benefit of relying only on convexity of the modified weights assigned to the jumps in rank in the rank sequences, not the exact formula. This could mean that the symmetry will hold in more general frameworks as well.
\end{remark}


\begin{thebibliography}{99}

\bibitem{topDim} J. Koncki, R. Rim\'anyi:
\textit{The main reasons for matrices multiplying to zero}, Linear Algebra Appl. \textbf{730}, pp. 587--602.,
https://doi.org/10.1016/j.laa.2025.11.001

\bibitem{symm} S. Pepin Lehalleur, R. Rim\'anyi:
\textit{Geometry of the fibers of the multiplication map of deep linear neural networks},
preprint at \url{https://rimanyi.web.unc.edu/wp-content/uploads/sites/9870/2024/11/Dimension_and\_RLCT_for_deep_linear_networks-1.pdf}

\end{thebibliography}
\end{document}